\documentclass[12pt]{amsart}

\newcommand{\newsection}[1]{\setcounter{equation}{0} \section{#1}}
\newcommand{\bea}{\begin{eqnarray}}
\newcommand{\eea}{\end{eqnarray}}

\newcommand{\cla}{\mathcal{A}}

\newcommand{\clh}{\mathcal{H}}
\newcommand{\clk}{\mathcal{K}}
\newcommand{\cll}{\mathcal{L}}
\newcommand{\clm}{\mathcal{M}}

\newcommand{\clq}{\mathcal{Q}}
\newcommand{\clr}{\mathcal{R}}

\newcommand{\clz}{\mathcal{Z}}

\newcommand{\raro}{\rightarrow}

\newcommand{\cl}[1]{\mathcal{#1}}

\def \qed {\hfill \vrule height6pt width 6pt depth 0pt}
\def\textmatrix#1&#2\\#3&#4\\{\bigl({#1 \atop #3}\ {#2 \atop #4}\bigr)}
\def\dispmatrix#1&#2\\#3&#4\\{\left({#1 \atop #3}\ {#2 \atop #4}\right)}
\newcommand{\be}{\begin{equation}}
\newcommand{\ee}{\end{equation}}
\newcommand{\ben}{\begin{eqnarray*}}
\newcommand{\een}{\end{eqnarray*}}

\newcommand{\NI}{\noindent}

\newcommand{\bi}{\begin{itemize}}
\newcommand{\ei}{\end{itemize}}

\def\5{{5\superprime}}

\newtheorem{Theorem}{\sc Theorem}
\newtheorem{Lemma}{\sc Lemma}
\newtheorem{Proposition}{\sc Proposition}
\newtheorem{Corollary}{\sc Corollary}
\newtheorem{Definition}{\sc Definition}
\newtheorem{Example}{\sc Example}
\newtheorem{Remark}{\sc Remark}
\newtheorem{Note}{\sc Note}
\newtheorem{Question}{\sc Question}
\newtheorem{ass}{\sc Assumption}
\newcommand{\bt}{\begin{Theorem}}
\def\beginlem{\begin{Lemma}}
\def\beginprop{\begin{Proposition}}
\def\begincor{\begin{Corollary}}
\def\begindef{\begin{Definition}}
\def\beginexamp{\begin{Example}}
\def\beginrem{\begin{Remark}}
\def\beginq{\begin{Question}}
\def\beginass{\begin{ass}}
\def\beginnote{\begin{Note}}
\newcommand{\et}{\end{Theorem}}
\def\endlem{\end{Lemma}}
\def\endprop{\end{Proposition}}
\def\endcor{\end{Corollary}}
\def\enddef{\end{Definition}}
\def\endexamp{\end{Example}}
\def\endrem{\end{Remark}}
\def\endq{\end{Question}}
\def\endass{\end{ass}}
\def\endnote{\end{Note}}

\begin{document}

\title{A note on semi-fredholm Hilbert modules}

\author[Douglas]{Ronald G. Douglas}
\author[Sarkar]{Jaydeb Sarkar}

\address{Texas A \& M University, College Station, Texas 77843, USA}

\email{rdouglas@math.tamu.edu, jsarkar@math.tamu.edu}

\thanks{This research was partially supported by a grant from the National Science Foundation.}

\keywords{Hilbert modules, quasi-free Hilbert modules, Fredholm tuple, Corona property}

\subjclass[2000]{47A13, 46E22, 46M20, 47B32}

\begin{abstract}
A classical problem in operator theory has been to determine the spectrum of Toeplitz-like operators on Hilbert spaces of vector-valued holomorphic functions on the open unit ball in $\mathbb{C}^m$. In this note we obtain necessary conditions for $k$-tuples of such operators to be Fredholm in the sense of Taylor and show they are sufficient in the case of the unit disk.
\end{abstract}

\maketitle

\newsection{Introduction}

A classical problem in operator theory is to determine the invertibility or the spectrum of Toeplitz operators on the Hardy space over the unit disk $\mathbb{D}$. Where the symbol or the defining function is continuous, the result is well known and due to Gohberg in the scalar case (see  \cite{G}) and Gohberg-Krein in the vector-valued case (see \cite{GK}). Generalizations of these results to other Hilbert spaces of holomorphic functions on the disk such as the Bergman space (see \cite{Ax}) or to the unit ball $\mathbb{B}^m$ (see \cite{Ven}) or other domains in $\mathbb{C}^m$ (see \cite{B}) have been studied during the past few decades. In the several variables context, the problem is not too interesting unless we start with a matrix-valued symbol or a $k$-tuple of operators and consider the Taylor spectrum or essential spectrum which involves the Koszul complex (see \cite{Tay}).

In this note we consider two problems, neither of which is new. However, we believe the results are more general and our methods provide a more constructive approach. Moreover, they  identify some questions in multi-variable operator theory (and algebra) indicating their importance in the spectral theory for $k$-tuples of vector-valued Toeplitz-like operators. Finally, the results suggest lines of investigation for generalizations of the classical Hilbert spaces of holomorphic functions.

\vspace{.5cm}

All the Hilbert spaces in this note are separable and are over
the complex field $\mathbb{C}$. For a Hilbert space $\clh$, we
denote the Banach space of all bounded linear operators by
$\cll(\clh)$.

We begin by recalling the definition of quasi-free Hilbert module over $A(\Omega)$
which was introduced in (\cite{DM1},\cite{DM2}) and which generalizes classical functional Hilbert space and is related to earlier ideas of Curto--Salinas \cite{C-S}.  Here $A(\Omega)$ is the uniform closure of functions holomorphic on a neighborhood of the closure of $\Omega$, a domain in $\mathbb{C}^m$. The Hilbert space $\clm$ is
said to be a bounded Hilbert module over $A(\Omega)$ if $\clm$ is a unital
module over $A(\Omega)$ with module map $A(\Omega) \times \clm \raro \clm$ such
that
$$\|\varphi f \|_{\cl M} \leq C \, \|\varphi\|_{A(\Omega)} \|f\|_{\cl M}$$
 for $\varphi$ in
$A(\Omega)$ and $f$ in $\clm$ and some $C \geq 1$. The Hilbert module is said to be contractive in case $C=1$. 

A Hilbert space $\clr$ is said to be a bounded quasi-free Hilbert module of rank $n$
over $A(\Omega)$, $1 \leq n \leq \infty$, if it is obtained as the completion of
the algebraic tensor product $A(\Omega) \otimes \ell^2_n$ relative to an inner
product such that:

(1) $\mbox{eval}_{\pmb{z}}\colon\ A(\Omega) \otimes l^2_n \raro l^2_n$ is bounded for $\pmb{z}$ in
$\Omega$ and locally uniformly bounded on $\Omega$;

(2) $\|\varphi(\sum \theta_i \otimes x_i)\|_{\clr} = \|\sum \varphi \theta_i \otimes
x_i\|_{\clr} \leq C \, \|\varphi\|_{A(\Omega)} \|\sum\theta_i \otimes x_i\|_{\clr}$ for
$\varphi$, $\{\theta_i\}$ in $A(\Omega)$ and $\{x_i\}$ in $\ell^2_n$ and some $C \geq 1$; and

(3) For $\{F_i\}$ a sequence in $A(\Omega) \otimes \ell^2_n$ which is Cauchy in
the $\clr$-norm, it follows that $\mbox{eval}_{\pmb{z}}(F_i) \raro 0$ for all $\pmb{z}$ in
$\Omega$ if and only if $\|F_i\|_{\clr} \raro 0$.

If $I_{\pmb{\omega}_0}$ denotes the maximal ideal of polynomials in $\mathbb{C}[\pmb{z}] = \mathbb{C}[z_1,\ldots, z_m]$ which  vanish at $\pmb{\omega}_0$ for some $\pmb{\omega}_0$ in $\Omega$, then the Hilbert module $\mathcal{M}$ is said to be semi-Fredholm at $\pmb{\omega}_0$ if $\dim \mathcal{M}/{I_{\pmb{\omega}_0}\cdot \mathcal{M}} = n$ is finite (cf. \cite{E}). In particular, note that $\mathcal{M}$ semi-Fredholm at $\pmb{\omega}_0$ implies that $I_{\pmb{\omega}_0}\mathcal{M}$ is a closed submodule of $\mathcal{M}$. Note that the notion of semi-Fredholm Hilbert module has been called regular by some authors.

One can show  that $\pmb{\omega}\to \mathcal{R}/I_{\pmb{\omega}}\cdot\mathcal{R}$ can be made into a rank $n$ Hermitian holomorphic vector bundle over $\Omega$ if $\mathcal{R}$ is semi-Fredholm at $\pmb{\omega}$ in $\Omega$, $\dim \mathcal{R}/{I_{\pmb{\omega}}\cdot \mathcal{R}}$ is constant $n$, and $\mathcal{R}$ is quasi-free, $1\le n<\infty$. Actually, all we need here is that the bundle obtained is real-analytic which is established in
(\cite{C-S}, Theorem 2.2).

A quasi-free Hilbert module of rank $n$ is a reproducing kernel Hilbert space
with the kernel
$$K(\pmb{w}, \pmb{z}) = \mbox{eval}_{\pmb{w}} \mbox{eval}_{\pmb{z}}^* \colon\ \Omega \times \Omega
\raro \cll(\ell^2_n).$$

\newsection{Necessary conditions}

Note that if $\clr$ is a bounded quasi-free Hilbert module over $A(\mathbb{B}^m)$ of finite multiplicity, then the module $\clr$ over $A(\mathbb{B}^m)$ extends to a bounded Hilbert module over $H^{\infty}(\mathbb{B}^m)$ (see Proposition 5.2 in \cite{DD}). Here $\mathbb{B}^m$ denotes the unit ball in $\mathbb{C}^m$. In particular, the multiplier space of $\clr$ is precisely $H^{\infty}(\mathbb{B}^m) \otimes \clm_n(\mathbb{C})$, since $\clr$ is the completion of $A(\Omega) \otimes_{\mbox{alg}} l^2_n$, by definition.

\begin{Proposition}\label{prop1}
Let $\clr$ be a contractive quasi-free Hilbert module over $A(\mathbb{B}^m)$ of finite multiplicity $n$ and $\{\varphi_1, \ldots, \varphi_k\}$ be a subset of $H^{\infty}(\mathbb{B}^m) \otimes \clm_n(\mathbb{C})$. If $(M_{\varphi_1}, \ldots, M_{\varphi_k})$ is a semi-Fredholm tuple, then there exists an $\epsilon >0$ and $1>\delta > 0$ such that $$\sum_{i=1}^{k} \varphi_i(\pmb{z}) \varphi_i(\pmb{z})^*  \geq \epsilon I_{\mathbb{C}^n} ,$$ for all $\pmb{z}$ satisfying $1>\|\pmb{z}\| \geq 1-\delta>0$. In particular, if the multiplicity of $\clr$ is one then $$\sum_{i=1}^{k} |\varphi_i(\pmb{z})|^2 \geq \epsilon,$$for all $\pmb{z}$ satisfying $1>\|\pmb{z}\| \geq 1-\delta$.
\end{Proposition}
\NI\textsf{Proof.} Let $K : \mathbb{B}^m \times \mathbb{B}^m \raro \clm_n(\mathbb{C})$ be the kernel function for the quasi-free Hilbert module $\clr$. By the assumption, the range of the row operator $M_{\Phi} = (M_{\varphi_1}, \ldots, M_{\varphi_k})$ in $\cll(\clr^k, \clr)$ has finite co-dimension; that is, $$\text{dim} [\clr/(M_{\varphi_1} \clr + \ldots + M_{\varphi_k} \clr)] < \infty,$$ and, in particular, $M_{\Phi}$ has closed range. Consequently, there is a finite rank projection $F$ such that $$M_{\Phi} M_{\Phi}^* + F = \sum_{i=1}^{k} M_{\varphi_i} M_{\varphi_i}^* + F : \clr \raro \clr$$ is bounded below. Therefore, there exists a $C>0$ such that $$< F K_{\pmb{z}}, K_{\pmb{z}}> + < \sum_{i=1}^{k} M_{\varphi_i} M_{\varphi_i}^* K_{\pmb{z}}, K_{\pmb{z}}> \; \geq C <K_{\pmb{z}},  K_{\pmb{z}}>,$$ for all $\pmb{z}$ in $\mathbb{B}^m$. Then $$  K_{\pmb{z}}^* \hat{F}(\pmb{z}) K_{\pmb{z}} + \sum_{i=1}^{k}  K_{\pmb{z}}^* M_{\varphi_i} M_{\varphi_i}^* K_{\pmb{z}} \geq C K_{\pmb{z}}^* K_{\pmb{z}},$$ and so $$ \hat{F}(z) I_{\mathbb{C}^n} + \sum_{i=1}^{k} \varphi_i(\pmb{z}) \varphi_i(\pmb{z})^* \geq C I_{\mathbb{C}^n},$$ for all $\pmb{z}$ in $\mathbb{B}^m$. Here $\hat{F} ({\pmb{z}})$ denotes the matrix-valued Berezin transform for the operator $F$ defined by $\hat{F} ({\pmb{z}}) = < F K_{\pmb{z}}|K_{\pmb{z}}|^{-1}, K_{\pmb{z}} |K_{\pmb{z}}|^{-1}>$ (see \cite{DD}, where the scalar case is discussed).  Using the known boundary behavior of the Berezin transform (see Theorem 3.2 in \cite{DD}), since $F$ is finite rank we have that $\|\hat{F}({\pmb{z}})\| \leq \frac{C}{2}$ for all $\pmb{z}$ such that $1 > \|\pmb{z}\| > 1- \delta$ for some $1 > \delta>0$ depending on $C$. Hence $$\sum_{i=1}^{k} \varphi_i(\pmb{z}) \varphi_i(\pmb{z})^* \geq \frac{C}{2},$$ for all $\pmb{z}$ such that $1>\|\pmb{z}\|>1-\delta >0$; which completes the proof. \qed

\vspace{0.1in}
A $k$-tuple of matrix-valued functions $(\varphi_1, \ldots, \varphi_k)$ in $H^{\infty}(\mathbb{B}^m) \otimes M_n(\mathbb{C})$ satisfying the conclusion of Proposition \ref{prop1} will be said to have the \textit{weak Corona property}.

In Theorem 8.2.6 in \cite{EP}, a version of Proposition \ref{prop1} is established in case $\clr$ is the Bergman module on $\mathbb{B}^m$.

The key step in this proof is the vanishing of the Berezin transform at the boundary of $\mathbb{B}^m$ for a compact operator. The proof of this statement depends on the fact that $K_{\pmb{z}} | K_{\pmb{z}}|^{-1}$ converges weakly to zero as $\pmb{z}$ approaches the boundary which rests on the fact that $\clr$ is contractive. This relation holds for many other domains such as ellipsoids $\Omega$ with the proof depending on the fact that the algebra $A(\Omega)$ is pointed in the sense of \cite{DD}.

It is an important question to decide if semi-Fredholm implies Fredholm in the context of Proposition \ref{prop1}. We will discuss this issue more at the end of the paper. However, the converse of this result is known (see Theorem 8.2.4 in \cite{EP} and pages 241-242) for the Bergman space for certain domains in $\mathbb{C}^m$. 

A necessary condition for the converse to hold for the situation in Proposition \ref{prop1} is for the essential spectrum of the $m$-tuple of co-ordinate multiplication operators to have essential spectrum equal to $\partial \mathbb{B}^m$, which is not automatic, but is true for the classical spaces.

\newsection{Sufficient conditions}

We will use the following fundamental result of Taylor (see \cite{Tay}, Lemma 1):

\begin{Lemma}\label{Tay}
Let $(T_1, \ldots, T_k)$ be in the center of an algebra $\cla$ contained in $\cll(\clh)$ such that there exists $(S_1, \ldots, S_k)$ in $\cla$  satisfying $\sum_{i=1}^{k} T_i S_i = I_{\clh}$. Then the Koszul complex for $(T_1, \ldots, T_k)$ is exact.
\end{Lemma}

Now we specialize to the case when $m=1$ where we can obtain a necessary and sufficient condition. Consider a contractive quasi-free Hilbert module $\clr$ over $A(\mathbb{D})$ of multiplicity one, which therefore has $H^{\infty}(\mathbb{D})$ as the multiplier algebra. It is well known that $H^{\infty}(\mathbb{D})$ satisfies the \textit{Corona property}; that is, a set $\{\varphi_1, \ldots, \varphi_k\}$ in $H^{\infty}(\mathbb{D})$ satisfies $\sum_{i=1}^{k} |\varphi_k(z)| \geq \epsilon$ for all $z$ in $\mathbb{D}$  for some $\epsilon >0$ if and only if there exist  $\{\psi_1, \ldots, \psi_k\} \subset H^{\infty}(\mathbb{D})$ such that $\sum_{i=1}^{k} \varphi \psi_i = 1$.

The following result is a complement to Proposition \ref{prop1}.
\begin{Proposition}
Let $\clr$ be a contractive quasi-free Hilbert module over $A(\mathbb{D})$ of multiplicity one and $\{\varphi_1, \ldots, \varphi_k\}$ be a subset of $H^{\infty}(\mathbb{D})$. Then the Koszul complex for the $k$-tuple $(M_{\varphi_1}, \ldots, M_{\varphi_k})$ on $\clr$ is exact if and only if $\{\varphi_1, \ldots, \varphi_k\}$ satisfies the Corona property.
\end{Proposition}
\NI \textsf{Proof.} If $\sum_{i=1}^{k} \varphi_i \psi_i = 1$ for some $\{\psi_1, \ldots, \psi_k\} \subset H^{\infty}(\mathbb{D})$, then the fact that $M_{\Phi}$ is Taylor invertible follows from Lemma \ref{Tay}. On the other hand, the last group of the Koszul complex is $\{0\}$ if and only if the row operator $M_{\varphi}$ in $\cll(\clr^k, \clr)$ is bounded below which, as before, shows that $\sum_{i=1}^{k} |\varphi_i(z)|$ is bounded below on $\mathbb{D}$. This completes the proof. \qed

\vspace{0.1in}

The missing step to extend the result from $\mathbb{D}$ to the open unit ball $\mathbb{B}^m$ is the fact that it is unknown if the Corona condition for $\{\varphi_1, \ldots, \varphi_k\}$ in $H^{\infty}(\mathbb{B}^m)$ is equivalent to the Corona property. Other authors have considered this kind of question (\cite{TW}) for the case of Hardy-like spaces for the polydisk and ball. See \cite{TW} for some recent results and references.

\begin{Theorem}\label{thm1}
Let $\clr$ be a contractive quasi-free Hilbert module over $A(\mathbb{D})$ of multiplicity one, which is semi-Fredholm at each point $z$ in $\mathbb{D}$. If $\{\varphi_1, \ldots, \varphi_k\}$ is a subset of $H^{\infty}(\mathbb{D})$, then the $k$-tuple $M_{\Phi} = (M_{\varphi_1}, \ldots, M_{\varphi_k})$ is semi-Fredholm if and only if it is Fredholm if and only if  $(\varphi_1, \ldots, \varphi_k)$ satisfies the weak Corona condition.
\end{Theorem}
\NI \textsf{Proof.} If $M_{\Phi}$ is semi-Fredholm, then by Proposition \ref{prop1} there exist $\epsilon >0$ and $1 > \delta>0$ such that $$\sum_{i=1}^{k} |\varphi_i(z)|^2 \geq \epsilon,$$ for all $z$ such that $1>|z|>1-\delta > 0$. Let $\clz$ be the set $$\clz = \{ z \; \text{in}\; \mathbb{D} : \varphi_i (z) = 0 \; \text{for all}\; i = 1, \ldots, k\}.$$ Since the functions $\{\varphi_i\}_{i=1}^{k}$ can not all vanish for $z$ satisfying $1>|z|>1-\delta$, it follows that  the cardinality of the set $\clz$ is finite and we assume that $\mbox{card} (\clz) = N$. Let $$\clz = \{ z_1, z_2, \ldots, z_N\}$$ and $l_j$ be the smallest order of the zero at $z_j$ for all $\varphi_j$ and $1 \leq j \leq k$. Let $B(z)$ be the finite Blaschke product with zero set precisely $\clz$ counting the multiplicities. If we define $\xi_i = \frac{\varphi_i}{B}$, then $\xi_i$ is in $H^{\infty}(\mathbb{D})$ for all $i=1, \ldots, k$. Since $\{\varphi_1, \ldots, \varphi_k\}$ satisfies the weak Corona property, we obtain $$\sum_{i=1}^{k} |\xi_i (z)|^2 \geq \epsilon$$ for all $z$ such that $1>|z|>1-\delta$. Note that $\{\xi_1, \ldots, \xi_n\}$ does not have any common zero and so $$\sum_{i=1}^{k} |\xi_i (z)|^2 \geq \epsilon$$ for all $z$ in $\mathbb{D}$. Therefore, $\{\xi_1, \ldots, \xi_k\}$ satisfies the Corona property and hence there exists $\{\eta_1, \ldots, \eta_k\}$, a subset of $H^{\infty}(\mathbb{D})$, such that $\sum_{i=1}^{k} \xi_i(z) \eta_i(z) = 1$ for all $z$ in $\mathbb{D}$. Thus, $\sum_{i=1}^{k} \varphi_i(z) \eta_i(z) = B$ for all $z$ in $\mathbb{D}$. This implies $\sum_{i=1}^{k} M_{\varphi_i} M_{\eta_i} = M_B$, and consequently, $$\sum_{i=1}^{k} \overline{M}_{\varphi_i} \overline{M}_{\eta_i} = \overline{M_B},$$ where $\overline{M}_{\varphi_i}$ is the image of $M_{\varphi_i}$ in the Calkin algebra, $\clq(\clr) = \cll(\clr)/\clk(\clr)$. But the assumption that $M_{z-w}$ is Fredholm for all $w$ in $\mathbb{D}$ yields that $M_B$ is Fredholm. Therefore,  $X = \sum_{i=1}^{k} \overline{M}_{\varphi_i} \overline{M}_{\eta_i}$ is invertible. Moreover, since $X$ commutes with the set $\{\overline{M}_{\varphi_1}, \ldots, \overline{M}_{\varphi_k}, \overline{M}_{\eta_1}, \ldots, \overline{M}_{\eta_k}\}$, it follows that $(M_{\varphi_1}, \ldots, M_{\varphi_k})$ is a Fredholm tuple, which completes the proof. \qed

\vspace{0.2in}

Although, the use of a finite Blaschke product allows one to preserve norms, a polynomial with the zeros of $\clz$ to the same multiplicity could be used. This would allow one to extend the Theorem to all domains in $\mathbb{C}$ for which the Corona theorem holds.

Our previous result extends to the case of finite multiplicity quasi-free Hilbert modules.

\begin{Theorem}\label{thm2}
Let $\clr$ be a contractive quasi-free Hilbert module over $A(\mathbb{D})$ of multiplicity $n$, which is semi-Fredholm at each point $z$ in $\mathbb{D}$ and let $\{\varphi_1, \ldots, \varphi_k\}$ be a subset of $H^{\infty}(\mathbb{D}) \otimes M_n(\mathbb{C})$. Then the $k$-tuple $M_{\Phi} = (M_{\varphi_1}, \ldots, M_{\varphi_k})$ is Fredholm if and only if it is semi-Fredholm if and only if $(\varphi_1, \ldots, \varphi_k)$ satisfies the weak Corona condition.
\end{Theorem}
\NI \textsf{Proof.} As before, the assumption that $M_{\Phi}$ is semi-Fredholm implies that there exists $\epsilon>0$ and $1>\delta>0$ such that
$$\sum_{i=1}^{k} \varphi_i(z) \varphi_i(z)^* \geq \epsilon I_{\mathbb{C}^n},$$ for all $z$ such that $1>\|z\|>1-\delta$. After taking the determinant, this inequality implies $$\sum_{i=1}^{k} | \mbox{det} \, \varphi_i(z)|^2 \geq \epsilon^n.$$ Using the same argument as in Theorem \ref{thm1}, we can find $\eta_1, \ldots, \eta_k$ in $H^{\infty}(\mathbb{D})$ and a finite Blaschke product $B$ such that $$\sum_{i=1}^{k} \eta_i(z) \, \mbox{det} \, \varphi_i(z) = B(z),$$ for all $z$ in $\mathbb{D}$. For $1 \leq i \leq k$, let $\hat{\varphi_i}(z)$ be the cofactor matrix function of $\varphi_i(z)$ which is used in Cramer's Theorem. Then$$\hat{\varphi_i}(z) \phi_i(z) = \phi_i(z) \hat{\varphi_i}(z) = \mbox{det}\, \varphi_i(z) \, I_{\mathbb{C}^n},$$ for all $z$ in $\mathbb{D}$ and $1 \leq i \leq k$. Note that this relation implies that the algebra generated by the set $\{M_{\varphi_1}, \ldots, M_{\varphi_k}, M_{\hat{\varphi_1}}, \ldots, M_{\hat{\varphi_k}} \}$ is commutative. Thus we obtain $$\sum_{i=1}^{k}  \phi_i(z) \, \eta_i(z) \, \hat{\varphi_i}(z) = B(z) I_{\mathbb{C}^n},$$ or $$\sum_{i=1}^{k} \phi_i(z) \hat{\eta_i}(z)  = B(z) I_{\mathbb{C}^n},$$ where $\hat{\eta_i}(z) = \eta_i(z) \hat{\varphi_i}(z)$ is in $H^{\infty}(\mathbb{D}) \otimes \clm_n(\mathbb{C})$ and $1 \leq i \leq k$. Therefore we have that $$\sum_{i=1}^{k} M_{\varphi_i} M_{\hat{\eta_i}}  = M_B,$$ and consequently, the proof follows immediately from the last part of the proof of Theorem \ref{thm1}. \qed

\vspace{0.1in}

\newsection{Further comments}

One reason we are able to obtain a converse in the one variable case is that we can represent the zero variety of the ideal generated by the functions in terms of a single function, the finite Blaschke product (or polynomial). This is not surprising since $\mathbb{C}[z]$ is a principal ideal domain. This is, of course, not true for $\mathbb{C}[z_1, \ldots, z_m]$ for $m>1$ and hence one would need (at least) a finite set of functions to determine the zero variety for the ideal generated by the functions. How to do that in an efficient manner and how to relate the Fredholmness of the $k$-tuple to that of this generating set is not clear but is the key to answering many such questions.

What is required involves two steps, both in the realm of algebra. The first we have already mentioned but the second is how to relate the generators to the Koszul complex. 

Let us consider one example of what might be possible. Consider the case in which the $p_1(\pmb{z}), \ldots, p_k(\pmb{z})$ are polynomials in $\mathbb{C}[z_1, z_2]$ so that $\pmb{0}$ is the only common zero. Assume that there are sets of polynomials $\{q_1(\pmb{z}), \ldots, q_k(\pmb{z})\}$ and $\{r_1(\pmb{z}), \ldots, r_k(\pmb{z})\}$ such that $$\sum_{i=1}^{k} p_i(\pmb{z}) q_i(\pmb{z}) = z_1^{k_1}$$ and $$\sum_{i=1}^{k} p_i(\pmb{z}) r_i(\pmb{z}) = z_2^{k_2},$$
for some positive integers $k_1$ and $k_2$.

Two questions now arise: 

\NI (1) Does the assumption that $(M_{p_1}, \ldots, M_{p_k})$ is semi-Fredholm with $\clz = \{\pmb{0}\}$ imply the existence of the subsets $\{r_1, \ldots, r_k\}$ and $\{q_1, \ldots, q_k\}$ of $\mathbb{C}[z_1, z_2]$? What if the functions $\{p_1, \ldots, p_k\}$ are in $H^{\infty}(\mathbb{B}^2)$ and we seek  $\{r_1, \ldots, r_k\}$ and $\{q_1, \ldots, q_k\}$ in $H^{\infty}(\mathbb{B}^2)$?

\NI (2) If the functions $\{r_1, \ldots, r_k\}$ and $\{q_1, \ldots, q_k\}$ exist and we assume that $(M_{z_1^{k_1}}, M_{z_2^{k_2}})$ acting on the quasi-free Hilbert module $\clr$ is Fredholm, does it follow that $(M_{p_1}, \ldots, M_{p_k})$ is also. 

These questions can be generalized to the case where one would need more than two polynomials to determine the zero variety, either because the dimension $m$ is greater than 2 or because $\clz$ contains more than one point. But answering these questions in the simple case discussed above would be good start.

After this note was written, J. Eschmeier informed the authors that both questions have an affirmative answer, at least when the zero variety is a single point.

\end{document}